%% file: main.tex
\documentclass[10.5pt,letterpaper,fleqn, twocolumn]{article}

\usepackage{graphicx}
\usepackage{makeidx}
\usepackage{multicol}
\usepackage{crop}
\usepackage{amstext,amsthm}
\usepackage{amsmath,amssymb}
\usepackage{bm}
\usepackage{natbib}

\usepackage{mhchem}

\usepackage{multicol}
\usepackage{pslatex}
\usepackage[strict]{changepage}
\newcommand{\be}{\begin{equation}}
\newcommand{\ee}{\end{equation}}
\newcommand{\bea}{\begin{eqnarray}}
\newcommand{\eea}{\end{eqnarray}}
\newcommand{\bne}{\begin{equation*}}
\newcommand{\ene}{\end{equation*}}
\newcommand{\bi}{\begin{itemize}}
\newcommand{\ei}{\end{itemize}}

\newcommand{\bbm}{\begin{bmatrix}}
\newcommand{\ebm}{\end{bmatrix}}

\newcommand{\mr}{\mathrm}

\newcommand{\diff}[2]{\frac{d#1}{d#2}}
\newcommand{\pdiff}[2]{\frac{\partial#1}{\partial#2}}
\newcommand{\stoch}[1]{\boldsymbol{#1}}

\usepackage[top=0.75in, bottom=0.5in, left=0.83in, right=0.83in]{geometry}
\usepackage[hang,splitrule]{footmisc}





\usepackage{titlesec}		
\usepackage{subcaption}		
\usepackage{indentfirst}	

\titleformat{\section}
{\normalfont\bfseries}
{}{0.0pt}{\parindent=0cm}[\parindent=0.5cm]

\titlespacing*{\section}
{0pt}{1.5\baselineskip}{0.5\baselineskip}

\titleformat{\subsection}
{\normalfont\itshape}
{}{0.0pt}{\parindent=0cm}[\parindent=0.5cm]

\titlespacing*{\subsection}
{0pt}{1.0\baselineskip}{0.5\baselineskip}

\begin{document}

\twocolumn[{%
	
\phantom\\
\vspace{0.5in}
\begin{center}
\Large{\textbf{State Estimation for Continuous-Discrete-Time Nonlinear Stochastic Systems}}\\
\end{center}
\vspace{0.2in}

\begin{center}
Marcus Krogh Nielsen$^{\mr{a,c}}$, Tobias K. S. Ritschel$^{\mr{a}}$, Ib Christensen$^{\mr{c}}$, Jess Dragheim$^{\mr{c}}$, \\ Jakob Kj{\o}bsted Huusom$^{\mr{b}}$, Krist V. Gernaey$^{\mr{b}}$, John Bagterp J{\o}rgensen$^{\mr{a,}}$\footnotemark

\vspace{0.10in}

$^{\mr{a}}$ Department of Applied Mathematics and Compute Science, DTU Compute, Technical University of Denmark, DK-2800 Kgs. Lyngby, Denmark.

\vspace{0.10in}
$^{\mr{b}}$ Department of Chemical and Biochemical Engineering, DTU Chemical Engineering, Technical University of Denmark, DK-2800 Kgs. Lyngby, Denmark.

\vspace{0.10in}
$^{\mr{c}}$ Unibio A/S, DK-4000 Roskilde, Denmark.
\end{center}

\begin{abstract}
	\input{parts/abstract}
\end{abstract}


\vspace{0.2in}
}]

\footnotetext[1]{\parindent=0cm \small{Corresponding author: J. B. J{\o}rgensen (E-mail: {\tt\small jbjo@dtu.dk}).}}

%
%


\input{parts/introduction}

\input{parts/modelStochasticDifferentialEquations}

\input{parts/methodNonlinearStateEstimation}
\input{parts/methodSubSecCDEKF}
\input{parts/methodSubSecCDUKF}
\input{parts/methodSubSecCDEnKF}
\input{parts/methodSubSecCDPF}
\input{parts/discussion}

\input{parts/exampleModifiedFourTankSystem}

\input{parts/conclusion}

%
%

\bibliographystyle{chicago}
\section{References}


\def\refname{}
\def\bibsection{}

\bibliography{references/bibliography.bib}%


%
%
%
%
%


\end{document}

%% file: parts/abstract.tex
State estimation incorporates the feedback in optimization based advanced process control systems and is very important for the performance of model predictive control. We describe the extended Kalman filter, the unscented Kalman filter, the ensemble Kalman filter, and a particle filter for continuous-discrete time nonlinear systems involving stochastic differential equations. Continuous-discrete time nonlinear systems is a natural way to model physical systems controlled by digital controllers. We implement the state estimation methods in Matlab, illustrate and evaluate their performance using simulations of the modified four-tank system. This system is non-stiff and the state estimation methods are implemented numerically using an explicit numerical integration scheme. We evaluate the accuracy of the state estimation methods in terms of the mean absolute percentage error over the simulation horizon. Each method successfully estimates the states and unmeasured disturbances of the simulated modified four-tank system. The key contribution is an overview and comparison of state estimation methods for continuous-discrete time nonlinear stochastic systems. This can guide efficient implementations.

%% file: parts/introduction.tex
\section{Introduction}
\label{sec:Introduction}

% ...
%
% ...
State estimation is widely applied in advanced process control (APC) systems, e.g. for monitoring, fault-detection, and as part of model predictive control (MPC). The objective of state estimation is to predict and reconstruct the states of a mathematical model using measurements from a physical system.
The Kalman filter provides optimal estimates for systems with Gaussian process and measurement noise, but is limited to system with linear dynamics \citep{kalman:1960a}.
For nonlinear systems, the exact evolution of the state distribution can be computed as the solution to the Fokker-Planck equation (Kolmogorov's forward equation). However, the Fokker-Planck equation is a partial differential equation where the number of dimensions equal the number of states in the system. The Fokker-Planck equation suffers from the curse of dimensionality and solving it is therefore impractical for systems with more than a few states \citep{jazwinski:2007}. This paper describes four approximate methods for state estimation, 1) the extended Kalman filter (EKF), 2) the unscented Kalman filter (UKF), 3) the ensemble Kalman filter (EnKF), and 4) a particle filter (PF).
% ...

% ...
%
% ...
In the EKF, the equations of the original Kalman filter are applied on a local linearisation of a nonlinear system \citep{rawlings:etal:2017}. The EKF is a computationally efficient method, but the quality of the estimates depend on the nonlinearity of the system \citep{frogerais:etal:2011}. Additionally, some stability issues may arise in relation to fixed step-size solutions 
\citep{bucy:joseph:2005,joergensen:etal:2007a}.
In the UKF, an unscented transformation is used as an approximation for the first two moments of the true nonlinear distribution. The unscented transformation is propagated through the nonlinear dynamics and each sigma-point is updated using observations from the physical system \citep{julier:uhlmann:2004}. For some systems, the UKF has shown higher accuracy than the EKF, while still being computationally efficient \citep{wan:vandermerwe:2000}. However, the UKF also suffers from inaccuracy in highly nonlinear systems. For the UKF, some of the issues pertaining to nonlinearity and numerical instability have been addressed in more recent contributions \citep{kandepu:etal:2008,devivo:etal:20217}.
In the EnKF, a set of particles, the ensemble, is randomly sampled from the state distribution and propagated through the nonlinear system dynamics. Each particle in the ensemble is updated separately using the Kalman filter update when a measurement becomes available. The state estimates are computed statistically from the ensembles \citep{gillijns:etal:2006,myrseth:omre:2010,roth:etal:2015}.
In PFs, a set of particles is sampled from the state distribution and propagated through the nonlinear system dynamics. When a measurement becomes available, the particles are resampled in accordance with their likelihood of being observed. The likelihoods are computed using the innovation posterior distribution. Similarly to the EnKF, the state estimates are determined statistically from the particles \citep{arulampalam:etal:2002,rawlings:bakshi:2006,tulsyan:etal:2016}.
The EKF and UKF provide efficient state estimation, but suffers from loss of accuracy for highly nonlinear systems. The EnKF and PFs provide a set of sampled particles from the true nonlinear distribution, but their computational efficiency depends on the number of particles required for the estimates to reach the desired accuracy. As a result of this, the EnKF and PF can be computationally inefficient for high dimensional systems.
% ...

% ...
%
% ...
The aim of this paper is to provide a condensed overview of available methods for state estimation in {\em continuous-discrete time} nonlinear stochastic systems that are described using stochastic differential equations (SDEs). The intention is to guide efficient implementation by providing an overview for the continuous-discrete nonlinear state estimation methods. 
The nonlinear state estimation methods can generally be separated into two steps; prediction and filtering. In the prediction step (time update), the system is propagated through time based on past information from a physical system. In the filtering step (measurement update), the state estimates are updated with the latest measurement information.
% ...

% ...
%
% ...
The paper is structured as follows. In Section \ref{sec:model}, we present the nonlinear continuous-discrete stochastic differential equation models used in simulation and state estimation. In Section \ref{sec:method}, we present the EKF, the UKF, the EnKF, and a PF, and finally present a discussion of the methods. In Section \ref{sec:example}, we present a numerical example of state estimation for a modified four-tank system. Finally, we present conclusions in Section \ref{sec:Conclusion}.

%% file: parts/modelStochasticDifferentialEquations.tex
\section{ Nonlinear Continuous-Discrete Stochastic Systems}
\label{sec:model}

We consider continuous-discrete systems, in which the system state is described by nonlinear continuous stochastic differential equations and measurements are taken at discrete points in time. The nonlinear continuous-discrete stochastic differential equation models are defined as
\begin{subequations} \label{eq::model}
	\begin{align}
		\begin{split}
			d \stoch{x}(t) 
				&= f( t, \stoch{x}(t), u(t), d(t), \theta ) dt
						\\
				&+ \sigma( t, \stoch{x}(t), u(t), d(t), \theta ) d \stoch{\omega}(t)
					,
		\end{split}
				\\
		\stoch{y}( t_{k} )
			&= h( t_{k}, \stoch{x}(t_{k}), \theta ) + \stoch{v}(t_{k})
				,
	\end{align}
\end{subequations}
where $f(\cdot)$ is the drift function, $\sigma(\cdot)$ is the diffusion function, and $h(\cdot)$ is the measurement function. The states are $\stoch{x}(t) \in \mathbb{R}^{ n_{x} }$, the inputs are $u(t) \in \mathbb{R}^{ n_{u} }$, the disturbances are $d(t) \in \mathbb{R}^{ n_{d} }$, the parameters are $\theta \in \mathbb{R}^{ n_{\theta} }$, and the measurements are $\stoch{y}(t_{k}) \in \mathbb{R}^{ n_{y} }$. The process noise $\stoch{\omega}(t) \in \mathbb{R}^{ n_{\omega}}$ is a standard Wiener process, such that $d \stoch{\omega}(t) \sim \mathcal{N} \left( 0, I dt \right)$, and $\stoch{v}(t_{k}) \sim \mathcal{N} \left( 0, R \right)$ is the measurement noise. The initial state is assumed to be distributed as $\stoch{x}_{0} \sim \mathcal{N} \left( \bar{\stoch{x}}_{0}, P_{0} \right)$.

%% file: parts/methodNonlinearStateEstimation.tex
\section{ State Estimation in Nonlinear Systems }
\label{sec:method}

In this section, we present methods for state estimation in continuous-discrete nonlinear systems \eqref{eq::model}. These estimators are called continuous-discrete estimators.

%% file: parts/methodSubSecCDEKF.tex
\subsection{ Continuous-discrete extended Kalman filter }
The EKF is initialised with the mean and covariance of the initial state of the system described in Section \ref{sec:model}
\begin{align}
	\hat{x}_{0|0}
		&= \bar{x}_{0}
			,
	&&P_{0|0}
		= P_{0}
			.
\end{align}
\paragraph{Time update:} In the time update, the mean and covariance are computed as the solution to the ordinary differential equations (ODEs) for $t \in \left[ t_{k}, t_{k+1} \right]$
\begin{subequations} \label{eq::ekfTimeUpdate}
	\begin{align} 
		\diff{ \hat{x}_{k}(t) }{ t }
			&= f( t, \hat{x}_{k}(t), u(t), d(t), \theta )
				,	\\
		\diff{ P_{k}(t) }{ t }
			&= A_{k}(t) P_{k}(t) + P_{k}(t) A_{k}^{T}(t) + \sigma_{k}(t) \sigma_{k}^{T}(t)
				,
	\end{align}
\end{subequations}
where $\hat{x}_{k}(t_{k}) = \hat{x}_{k|k}$, and $P_{k}(t_{k}) = P_{k|k}$. $A_{k}(t) = \pdiff{ f }{ x }( t, \hat{x}_{k}(t), u(t), d(t), \theta )$ and $\sigma_{k}(t) = \sigma( t, \hat{x}_{k}(t), u(t), d(t), \theta )$. Alternatively, the covariance update can be represented and solved on integral form as presented by \cite{joergensen:etal:2007b}. The mean and covariance estimates are
\begin{align}
	\hat{x}_{k+1|k}
		&= \hat{x}_{k}( t_{k+1} )
			,
	&&P_{k+1|k}
		= P_{k}( t_{k+1} )
			.
\end{align}
\paragraph{Measurement update:} In the measurement update, we compute the innovation and its covariance as
\begin{align}
	e_{k} 
		&= y_{k} - \hat{y}_{k|k-1}
			,
	&&R_{e,k}
		= C_{k} P_{k|k-1} C_{k}^{T} + R
			,
\end{align}
where
\begin{align}
	\hat{y}_{k|k-1} 
		&= h( t_{k}, \hat{x}_{k|k-1}, \theta )
			,
	&&C_{k} 
		= \pdiff{ h }{ x }( t_{k}, \hat{x}_{k|k-1}, \theta )
			.
\end{align}
The Kalman gain is computed as
\begin{align}
	\begin{split}
		K_{f_{x},k}
			&= P_{k|k-1} C_{k}^{T} R_{e,k}^{-1}
				.
	\end{split}
\end{align}
The mean and covariance estimates are computed as
\begin{subequations}
	\begin{align}
		\hat{x}_{k|k}
			&= \hat{x}_{k|k-1} + K_{f_{x},k} e_{k}
				,	\label{eq::xFilt}	\\
		P_{k|k}
			&= P_{k|k-1} - K_{f_{x},k} R_{e,k} K_{f_{x},k}^{T}
					\label{eq::PFilt}	\\
			&= \left( I - K_{f_{x},k} C_{k} \right) P_{k|k-1} \left( I - K_{f_{x},k} C_{k} \right)^{T} + K_{f_{x},k} R K_{f_{x},k}^{T}
				, \label{eq::PFiltJoseph}
	\end{align}
\end{subequations}
where (\ref{eq::PFiltJoseph}), Joseph's stabilising form, is numerically stable.

%% file: parts/methodSubSecCDUKF.tex
\subsection{ Continuous-discrete unscented Kalman filter }
The unscented Kalman filter is initialised with the mean and covariance of the initial state of the system described in Section \ref{sec:model}
\begin{align}
	\hat{x}_{0|0}
		&= \bar{x}_{0}
			,
	&&P_{0|0}
		= P_{0}
			.
\end{align}

\paragraph{Time update:} In the time update, we compute the parameters
\begin{align} \label{eq::ukfParsTimeUpdate}
	\bar{c}
		&= \alpha^{2} \left( \bar{n} + \kappa \right)
			,	\\
	\bar{\lambda}
		&= \alpha^{2} \left( \bar{n} + \kappa \right) - \bar{n}
			,
\end{align}
where $\alpha \in ] 0, 1 ]$, $\kappa \in \left[ 0, \infty \right[$, and $\bar{n} = n_{x} + n_{\omega}$. We compute the sigma-point weights
\begin{subequations}
	\begin{align}
		\bar{W}_{m}^{(0)}
			&= \frac{ \bar{\lambda} }{ \bar{n} + \bar{\lambda} }
				,	\\
		\bar{W}_{c}^{(0)}
			&= \frac{ \bar{\lambda} }{ \bar{n} + \bar{\lambda} } + 1 - \alpha^{2} + \beta
				,	\\
		\bar{W}_{m}^{(i)}
			&= \bar{W}_{c}^{(i)}
			= \frac{ 1 }{ 2 \left( \bar{n} + \bar{\lambda} \right) }
				,
	\end{align}
\end{subequations}
for $i \in \{ 1, 2, \dots, 2 \bar{n} \}$ and where $\beta \in \left[ 0, \infty \right[$ ($\beta = 2$ optimal for Gaussian distributions). We sample deterministically a set of $2 \bar{n} + 1$ sigma-points. For propagation through the deterministic dynamics (ODE), we compute $2 n_{x} + 1$ sigma-points
\begin{subequations}
	\begin{align}
		\hat{x}_{k|k}^{(0)}
			&= \hat{x}_{k|k}
				,	\\
		\hat{x}_{k|k}^{(i)}
			&= \hat{x}_{k|k} + \sqrt{ \bar{c} } \left( \sqrt{ P_{k|k} } \right)_{i}
				,	\\
		\hat{x}_{k|k}^{(n_{x} + i)}
			&= \hat{x}_{k|k} - \sqrt{ \bar{c} } \left( \sqrt{ P_{k|k} } \right)_{i}
				,
	\end{align}
\end{subequations}
for $i \in \{ 1, 2, \dots, n_{x} \}$. $\left( \sqrt{P_{k|k}} \right)_{i}$ denotes the $i$'th column of the Cholesky decomposition of the covariance. For propagation through the stochastic dynamics, we compute $2n_{\omega}$ sigma-points
\begin{align}
	\hat{x}_{k|k}^{(2n_{x} + i)}
		&= \hat{x}_{k|k}
			,
\end{align}
for $i \in \{ 1, 2, \dots, 2n_{\omega} \}$. Additionally, we compute the process noise
\begin{subequations}
    \begin{align}
    	d \omega_{k}^{(2n_{x} + i)}(t)
    		&= \sqrt{ \bar{c} ~ dt } \left( I \right)_{i}
    			,	\\
    	d \omega_{k}^{(2n_{x} + n_{\omega} + i)}(t)
    		&= - \sqrt{ \bar{c} ~ dt } \left( I \right)_{i}
    			,
    \end{align}
\end{subequations}
where $i \in \{ 1, 2, \dots, n_{\omega} \}$. We propagate the first sigma-points through the deterministic dynamics for $t \in \left[ t_{k}, t_{k+1} \right]$ and compute the predictions as the solution to
\begin{align} 	\label{eq::UKFPredODE}
	d \hat{x}_{k}^{(i)}(t)
		&= f( t, \hat{x}_{k}^{(i)}(t), u(t), d(t), \theta ) dt
			,
\end{align}
for $\hat{x}_{k}^{(i)}(t_{k}) = \hat{x}_{k|k}^{(i)}$ and $i \in \{ 0, 1, \dots, 2n_{x} \}$. We similarly propagate the remaining sigma-points through the stochastic dynamics for $t \in \left[ t_{k}, t_{k+1} \right]$ and compute the predictions as the solution to
\begin{align} \label{eq::UKFPredSDE}
	\begin{split}
		d \hat{x}_{k}^{(i)}(t)
			&= f( t, \hat{x}_{k}^{(i)}(t), u(t), d(t), \theta ) dt
					\\
			&+ \sigma( t, \hat{x}_{k}^{(i)}(t), u(t), d(t), \theta ) d \omega_{k}^{(i)}(t)
				,
	\end{split}
\end{align}
where $\hat{x}_{k}^{(i)}(t_{k}) = \hat{x}_{k|k}^{(i)}$ and $i \in \{ 2n_{x} + 1, 2n_{x} + 2, \dots, 2n_{x} + 2 n_{\omega} \}$. The predictions are computed as the solution to \eqref{eq::UKFPredODE} and \eqref{eq::UKFPredSDE}, as $\hat{x}_{k+1|k}^{(i)} = \hat{x}_{k}^{(i)}(t_{k+1})$. The mean and covariance estimates are computed as
\begin{subequations}
	\begin{align}
		\hat{x}_{k+1|k}
			&= \sum_{ i = 0 }^{ 2 \bar{n} } \bar{W}_{m}^{(i)} \hat{x}_{k+1|k}^{(i)}
				,	\\
		P_{k+1|k}
			&= \sum_{ i = 0 }^{ 2 \bar{n} } \bar{W}_{c}^{(i)} \left( 	\hat{x}_{k+1|k}^{(i)} - \hat{x}_{k+1|k} \right)\left( \hat{x}_{k+1|k}^{(i)} - \hat{x}_{k+1|k} \right)^{T}
				.
	\end{align}
\end{subequations}

\paragraph{Measurement update:} In the measurement update, we compute the parameters
\begin{align}
	c
		&= \alpha^{2} \left( n_{x} + \kappa \right)
			,
	&&\lambda
		= \alpha^{2} \left( n_{x} + \kappa \right) - n_{x}
			.
\end{align}
We compute the sigma-point weights
\begin{subequations}
	\begin{align}
		W_{m}^{(0)}
			&= \frac{ \lambda }{ n_{x} + \lambda }
				,	\\
		W_{c}^{(0)}
			&= \frac{ \lambda }{ n_{x} + \lambda } + 1 - \alpha^{2} + \beta
				,	\\
		W_{m}^{(i)}
			&= W_{c}^{(i)}
			= \frac{ 1 }{ 2 \left( n_{x} + \lambda \right) }
				,
	\end{align}
\end{subequations}
for $i \in \{ 1, 2, \dots, 2 n_{x} \}$. We compute a set of $2n_{x} + 1$ deterministically sampled sigma-points
\begin{subequations}
	\begin{align}
		\hat{x}_{k|k-1}^{(0)}
			&= \hat{x}_{k|k-1}
				,	\\
		\hat{x}_{k|k-1}^{(i)}
			&= \hat{x}_{k|k-1} + \sqrt{c} \left( \sqrt{ P_{k|k-1} } \right)_{i}
				,	\\
		\hat{x}_{k|k-1}^{(n_{x}+i)}
			&= \hat{x}_{k|k-1} - \sqrt{c} \left( \sqrt{ P_{k|k-1} } \right)_{i}
				,
	\end{align}
\end{subequations}
for $i \in \{ 1, 2, \dots, n_{x} \}$. We compute the innovation as
\begin{align}
	e_{k}
		&= y_{k} - \hat{y}_{k|k-1}
			,
\end{align}
where the prediction of the measurement prediction is computed as
\begin{align}
	\hat{y}_{k|k-1}
		&= \hat{z}_{k|k-1}
		= \sum_{ i = 0 }^{ 2 n_{x} } W_{m}^{(i)} \hat{z}_{k|k-1}^{(i)}
			,
\end{align}
for $\hat{z}_{k|k-1}^{(i)} = h( t_{k}, \hat{x}_{k|k-1}^{(i)}, \theta )$. We compute the covariance and cross-covariance information from the sigma-points
\begin{small}
	\begin{subequations}
		\begin{align}
			R_{zz,k|k-1}
				&= \sum_{ i = 0 }^{ 2 n_{x} } W_{c}^{(i)} \left( \hat{z}_{k|k-1}^{(i)} - \hat{z}_{k|k-1} \right) \left( \hat{z}_{k|k-1}^{(i)} - \hat{z}_{k|k-1} \right)^{T}
					,	\\
			R_{e,k}
				&= R_{yy,k|k-1}
				= R_{zz,k|k-1} + R
					,	\\
			R_{xy,k|k-1}
				&= \sum_{ i = 0 }^{ 2 n_{x} } W_{c}^{(i)} \left( \hat{x}_{k|k-1}^{(i)} - \hat{x}_{k|k-1} \right) \left( \hat{z}_{k|k-1}^{(i)} - \hat{z}_{k|k-1} \right)^{T}
					.
		\end{align}
	\end{subequations}
\end{small}
The Kalman gain is computed as
\begin{align}
		K_{f_{x},k}
			&= R_{xy,k|k-1} R_{e,k}^{-1}
				.
\end{align}
The mean and covariance estimates are computed as
\begin{subequations}
	\begin{align}
		\hat{x}_{k|k}
			&= \hat{x}_{k|k-1} + K_{f_{x},k} e_{k}
				,	\\
		P_{k|k}
			&= P_{k|k-1} - K_{f_{x},k} R_{e,k} K_{f_{x},k}^{T}
				.
	\end{align}
\end{subequations}

%% file: parts/methodSubSecCDEnKF.tex
\subsection{ Continuous-discrete ensemble Kalman filter }
The ensemble Kalman filter is initialised with a set of particles, the ensemble, sampled from the initial state distribution from \eqref{sec:model}. The initial state ensemble is denoted $\{ \hat{x}_{0|0}^{(i)} \}_{ i = 1 }^{ N_{p} }$.

\paragraph{Time update:} In the time update, each particle in the ensemble is propagated through the system dynamics. The prediction ensemble is computed as the solution to
\begin{align}
	\begin{split}
		d \stoch{x}_{k}^{(i)}(t)
			&= f( t, \stoch{x}_{k}^{(i)}(t), u(t), d(t), \theta ) dt
					\\
			&+ \sigma( t, \stoch{x}_{k}^{(i)}(t), u(t), d(t), \theta ) d \stoch{\omega}_{k}(t)
				,
	\end{split}
\end{align}
for $i \in \{ 1, 2, \dots, N_{p} \}$ and $t \in \left[ t_{k}, t_{k+1} \right]$. The initial value is $x_{k}^{(i)} = \hat{x}_{k|k}^{(i)}$. The set of solutions, $\hat{x}_{k+1|k}^{(i)} = x_{k}^{(i)}(t_{k+1})$, gives rise to the prediction ensemble $\{ \hat{x}_{k+1|k}^{(i)} \}_{ i = 1 }^{ N_{p} }$. The mean and covariance estimates are computed as
\begin{small}
	\begin{subequations}
		\begin{align}
			\hat{x}_{k+1|k}
				&= \frac{ 1 }{ N_{p} } \sum_{ i = 1 }^{ N_{p} } \hat{x}_{k+1|k}^{(i)}
					,	\\
			P_{k+1|k}
				&= \frac{ 1 }{ N_{p} - 1 } \sum_{ i = 1 }^{ N_{p} } \left( \hat{x}_{k+1|k}^{(i)} - \hat{x}_{k+1|k} \right) \left( \hat{x}_{k+1|k}^{(i)} - \hat{x}_{k+1|k} \right)^{T}
					.
		\end{align}
	\end{subequations}
\end{small}

\paragraph{Measurement update:} In the measurement update, we compute the ensemble of predictions, $\{ \hat{z}_{k|k-1}^{(i)} \}_{ i = 1 }^{ N_{p} }$, where $z_{k|k-1}^{(i)} = h( t_{k}, \hat{x}_{k|k-1}^{(i)}, \theta )$, for $i \in \{ 1, 2, \dots, N_{p} \}$. Furthermore, we compute the mean and covariance of the measurement distribution and cross-covariance of states and measurements, as
\begin{footnotesize}
	\begin{subequations}
		\begin{align}
			\hat{y}_{k|k-1}
				&= \hat{z}_{k|k-1}
				= \frac{ 1 }{ N_{p} } \sum_{ i = 1 }^{ N_{p} } \hat{z}_{k|k-1}^{(i)}
					,	\\
			R_{zz,k|k-1}
				&= \frac{ 1 }{ N_{p} - 1 } \sum_{ i = 1 }^{ N_{p} } \left( \hat{z}_{k|k-1}^{(i)} - \hat{z}_{k|k-1} \right) \left( \hat{z}_{k|k-1}^{(i)} - \hat{z}_{k|k-1} \right)^{T}
					,	\\
			R_{yy,k|k-1}
				&= R_{zz,k|k-1} + R
					,	\\
			R_{xy,k|k-1}
				&= \frac{ 1 }{ N_{p} - 1 } \sum_{ i = 1 }^{ N_{p} } \left( \hat{x}_{k|k-1}^{(i)} - \hat{x}_{k|k-1} \right) \left( \hat{y}_{k|k-1}^{(i)} - \hat{y}_{k|k-1} \right)^{T}
					,
		\end{align}
	\end{subequations}
\end{footnotesize}
and we compute samples from measurement distribution, as
\begin{align}
	y_{k}^{(i)}
		&= y_{k} + v_{k}^{(i)}
			,
\end{align}
where $v_{k}^{(i)}$ are realisations of the measurement noise, $\stoch{v}_{k} \sim \mathcal{N}( 0, R )$. The innovations are computed for each particle in the measurement ensemble, as
\begin{align}
	e_{k}^{(i)}
		&= y_{k}^{(i)} - \hat{z}_{k|k-1}^{(i)}
			.
\end{align}
The Kalman gain is computed as
\begin{align}
	K_{f_{x},k}
		&= R_{xy,k|k-1} R_{yy,k|k-1}^{-1}
			.
\end{align}
The filtered state ensemble, $\{ \hat{x}_{k|k}^{(i)} \}_{ i = 1 }^{ N_{p} }$, is computed as
\begin{align}
	\hat{x}_{k|k}^{(i)}
		&= \hat{x}_{k|k-1}^{(i)} + K_{f_{x},k} e_{k}^{(i)}
			.
\end{align}
The mean and covariance estimates are computed as
\begin{subequations}
	\begin{align}
		\hat{x}_{k|k}
			&= \frac{ 1 }{ N_{p} } \sum_{ i = 1 }^{ N_{p} } \hat{x}_{k|k}^{(i)}
				,	\\
		P_{k|k}
			&= \frac{ 1 }{ N_{p} - 1 } \sum_{ i = 1 }^{ N_{p} } \left( \hat{x}_{k|k}^{(i)} - \hat{x}_{k|k} \right) \left( \hat{x}_{k|k}^{(i)} - \hat{x}_{k|k} \right)^{T}
				.
	\end{align}
\end{subequations}

%% file: parts/methodSubSecCDPF.tex
\subsection{ Continous-discrete particle filter }
The particle filter is initialised with a set of particles sampled from the initial state distribution from \eqref{sec:model}. The initial set of particles is denoted $\{ \hat{x}_{0|0}^{(i)} \}_{ i = 1 }^{ N_{p} }$.

\paragraph{Time update:} In the time update, each particle is propagated through the nonlinear system dynamics. The set of predicted particles is computed as the solution to
\begin{align}
	\begin{split}
		d \stoch{x}_{k}^{(i)}(t)
			&= f( t, \stoch{x}_{k}^{(i)}(t), u(t), d(t), \theta ) dt
					\\
			&+ \sigma( t, \stoch{x}_{k}^{(i)}(t), u(t), d(t), \theta ) d \stoch{\omega}_{k}(t)
				,
	\end{split}
\end{align}
for $i \in \{ 1, 2, \dots, N_{p} \}$ and $t \in \left[ t_{k}, t_{k+1} \right]$. The initial value $x_{k}^{(i)} = \hat{x}_{k|k}^{(i)}$. The set of solutions, $\hat{x}_{k+1|k}^{(i)} = x_{k}^{(i)}(t_{k+1})$, gives rise to the prediction set $\{ \hat{x}_{k+1|k}^{(i)} \}_{ i = 1 }^{ N_{p} }$. The mean and covariance estimates are computed as
\begin{small}
	\begin{subequations}
		\begin{align}
			\hat{x}_{k+1|k}
				&= \frac{ 1 }{ N_{p} } \sum_{ i = 1 }^{ N_{p} } \hat{x}_{k+1|k}^{(i)}
					,	\\
			P_{k+1|k}
				&= \frac{ 1 }{ N_{p} - 1 } \sum_{ i = 1 }^{ N_{p} } \left( \hat{x}_{k+1|k}^{(i)} - \hat{x}_{k+1|k} \right) \left( \hat{x}_{k+1|k}^{(i)} - \hat{x}_{k+1|k} \right)^{T}
					.
		\end{align}
	\end{subequations}
\end{small}

\paragraph{Measurement update:} In the measurement update, we compute the set of measurement predictions, $\{ \hat{z}_{k|k-1}^{(i)} \}_{ i = 1 }^{ N_{p} }$, where $\hat{z}_{k|k-1}^{(i)} = h( t_{k}, \hat{x}_{k|k-1}^{(i)}, \theta )$, for $i \in \{ 1, 2, \dots, N_{p} \}$. The innovations are computed for each particle, as
\begin{align}
	e_{k}^{(i)}
		&= y_{k} - \hat{z}_{k|k-1}^{(i)}
			,
\end{align}
for $i \in \{ 1, 2, \dots, N_{p} \}$. We compute a set of likelihood weights for each particle, arising from the posterior distribution of the innovations
\begin{align}
	\tilde{w}_{k}^{(i)}
		&= \frac{ 1 }{ \sqrt{ 2 \pi^{ n_{y} } \left| R \right| } } 
			\exp \left( - \frac{1}{2} \left( e_{k}^{(i)} \right)^{T} R^{-1} e_{k}^{(i)} \right)
				,
\end{align}
where $|R|$ denotes the determinant of $R$, and normalise
\begin{align}
	w_{k}^{(i)}
		&= \frac{ \tilde{w}_{k}^{(i)} }{ \sum_{ j = 1 }^{ N_{p} } \tilde{w}_{k}^{(j)} }
			,
\end{align}
for $i \in \{ 1, 2, \dots, N_{p} \}$. The set of particles are then resampled in accordance with their likelihood respective weights. For a single realisation of a uniform distribution, $q_{1} \sim \mathcal{U} \left[ 0, 1 \right]$, we compute a set of ordered resampling points
\begin{align}
	q_{k}^{(i)}
		&= \frac{ (i-1) + q_{1} }{ N_{p} }
			,
\end{align}
for $i \in \{ 1, 2, \dots, N_{p} \}$. We resample the particles by storing $m^{(i)}$ copies of each particle, $\hat{x}_{k|k-1}^{(i)}$, in the set. The indicies for the resampled particles, $l$, are chosen such that $q_{k}^{(l)} \in \left] s^{(i-1)}, s^{(i)} \right]$, where $s^{(i)} = \sum_{ j = 1 }^{ i } w_{k}^{(j)}$. Particles with relatively high likelihood may appear several times in the resampled set and particles with relatively low likelihood may not appear at all. The resampled set is denoted as $\{ \hat{x}_{k|k}^{(i)} \}_{ i = 1 }^{ N_{p} }$. The mean and covariance estimates are computed as
\begin{subequations}
	\begin{align}
		\hat{x}_{k|k}
			&= \frac{ 1 }{ N_{p} } \sum_{ i = 1 }^{ N_{p} } \hat{x}_{k|k}^{(i)}
				,	\\
		P_{k|k}
			&= \frac{ 1 }{ N_{p} - 1 } \sum_{ i = 1 }^{ N_{p} } \left( \hat{x}_{k|k}^{(i)} - \hat{x}_{k|k} \right) \left( \hat{x}_{k|k}^{(i)} - \hat{x}_{k|k} \right)^{T}
				.
	\end{align}
\end{subequations}

%% file: parts/discussion.tex
\subsection{Discussion of methods}
\label{subsec:discussion}
The EKF is a computationally efficient method for systems with a moderate number of states, while it is infeasible for systems with a very large number of states. The complexity in implementing the method is largely determined by computational aspects of solving the initial value problem \eqref{eq::ekfTimeUpdate} and issues related to computation of the covariance matrix. The accuracy of the EKF depends on how well the assumption of local linearity holds. This means that for highly nonlinear systems with relatively long sampling intervals, the EKF may perform poorly, as the assumptions pertaining to the propagation of the expectation and covariance will not hold.
The UKF is comparable to the EKF in terms of its computational requirements. The computational efficiency partly arises by utilising the unscented transformation, where the number of deterministically sampled particles scales linearly with the state dimension, instead of randomly sampling particles, as is the case for other particle filters. The time update of the UKF is simple to implement, as it simply involves propagating a set of particles forward in time. For linear Gaussian systems, the UKF and EKF provide equivalent solutions. However, for nonlinear systems, the UKF propagates the particles through the true nonlinear system dynamics and therefore may capture more information than the EKF. The particle filters, i.e. the EnKF and PF, have computational efficiency which depends on the tuning, i.e. the size of the sample set. They suffer from the curse of dimensionality, as the sampling size required increases with the state dimension. However, the predictions more closely resemble the true nonlinear distribution as the sampling size increases, at the cost of computational efficiency. This means, that for highly nonlinear systems the EnKF and PF may capture more information than the EKF and UKF, but at the cost of computational efficiency. Nevertheless, the EnKF is often used for large-scale systems, but with few samples.

%% file: parts/exampleModifiedFourTankSystem.tex
\section{ Example -- Modified Four-Tank System (MFTS)}
\label{sec:example}

The modified four tank system is modelled by a set of ODEs describing mass balances as presented by \cite{azam:joergensen:2015}. The model is further modified by modelling the stochastic disturbances explicitly as states. The disturbances are governed by the stochastic processes
\begin{subequations} \label{eq::disturbances}
	\begin{align}
		d \stoch{F}_{3}(t)
			&= \lambda_{1} \left( \bar{F}_{3}(t) - \stoch{F}_{3}(t) \right) dt + \sigma_{1} d \stoch{\omega}_{1}(t)
				,	\\
		d \stoch{F}_{4}(t)
			&= \lambda_{2} \left( \bar{F}_{4}(t) - \stoch{F}_{4}(t) \right) dt + \sigma_{2} d \stoch{\omega}_{2}(t)
				.
	\end{align}
\end{subequations}
The resulting system is described by a continuous-discrete nonlinear system as described in \eqref{eq::model}.The performance of each state estimation method is evaluated in terms of the mean absolute percentage error (MAPE) , such that
\begin{align}
	MAPE
		&= \frac{ 1 }{ n N }\sum_{ k = 1 }^{ N }  \sum_{ i = 1 }^{ n } \left| \frac{ x_{i,k} - \hat{x}_{i,k} }{ x_{i,k} } \right|
		    ,
\end{align}
where $N$ is the number of observations and $n$ is the dimension of the state. The MAPE is computed separately for the states representing the liquid mass and the state representing the disturbances, as MAPE$_{x}$ and MAPE$_{d}$ respectively.

\subsection{Simulation example}
\begin{figure*}[tb]
	\centering
	\includegraphics[width=\textwidth]{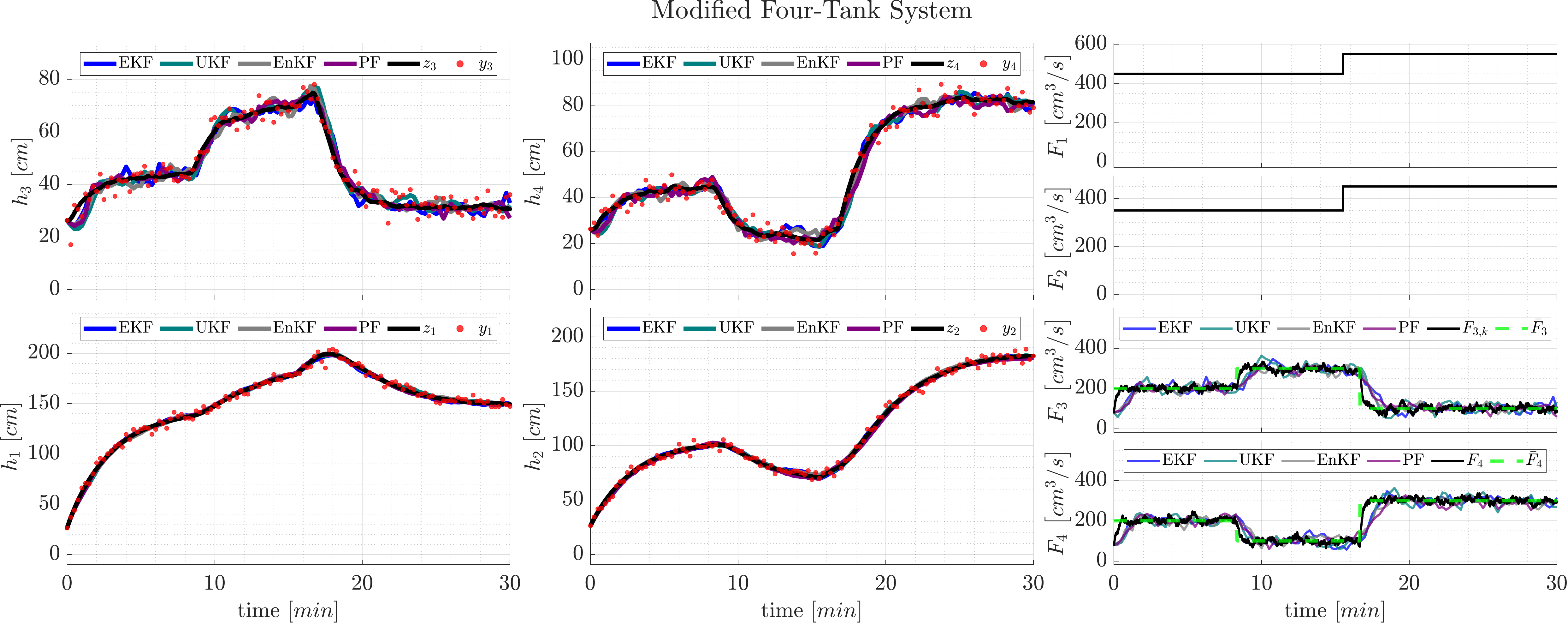}
	\caption{State estimation for simulation of the modified four-tank system. For the simulation model $\bar{F}_{i} \in \{ 100, 200, 300 \}$ for $i \in \{ 3, 4 \}$. The nominal values are kept constant at $150$ for the state estimators.}
	\label{fig::results}
\end{figure*}

Fig. \ref{fig::results} illustrates the simulation of the modified four tank system and Table \ref{tab::results} describes the results of the example. We simulate for $30$ minutes with $120$ equidistant samples. The simulation and estimation are computed with $1000$ and $100$ equidistant steps between samples, respectively. The UKF has the parameter set $\left[ \beta, \alpha, \kappa \right] = \left[ 2.0, 0.001, 0.0 \right]$. The EnKF and PF has particle set sizes of $250$ and $1000$, respectively. The disturbances are modelled with $\lambda_{1} = \lambda_{2} = 0.1$ and $\sigma_{1} = \sigma_{2} = 5.0$ for the simulation. For the EKF, EnKF, and PF, $\sigma_{1} = \sigma_{2} = 5.0$ and for the UKF $\sigma_{1} = \sigma_{2} = 1.0$. $\lambda_{1} = \lambda_{2} = 0.0$ for the EKF and UKF and $\lambda_{1} = \lambda_{2} =$ 2.0e-3 for the EnKF and PF. From the results presented in Fig. \ref{fig::results} and Table \ref{tab::results}, we see many of the properties described in the discussion of Section \ref{sec:method}. The EKF and UKF are demonstrated the be the most computationally efficient methods, where EKF seem to be outperforming UKF in this particular numerical experiment. Furthermore, the EnKF and PF show better accuracy both in estimating state and disturbance variables, but at the cost of lower computational efficiency.

\begin{table}[tb]
	\centering
	\caption{run-times for time update (TU) and measurement update (MU), and MAPE for states (MAPE$_{x}$) and disturbances (MAPE$_{d}$).}
	\label{tab::results}
	\begin{footnotesize}
	\begin{tabular}{r | l  l  l  l}
		\bf{name}         & \bf{EKF}        & \bf{UKF}        & \bf{EnKF}        & \bf{PF}      \\ \hline
        time TU [s]       & 3.09e-01        & 2.90e+00        & 3.38e+01         & 1.36e+02     \\     
        time MU [s]       & 1.22e-02        & 4.14e-02        & 2.30e-01         & 1.05e+00     \\
        MAPE$_{x}$ [\%]   & 2.55e+00        & 2.97e+00        & 2.35e+00         & 2.40e+00     \\
        MAPE$_{d}$ [\%]   & 1.57e+01        & 1.75e+01        & 1.47e+01         & 1.37e+01 
	\end{tabular}
	\end{footnotesize}
\end{table}

%% file: parts/conclusion.tex
\section{Conclusion}
\label{sec:Conclusion}
We present four methods for state estimation in continuous-discrete nonlinear systems involving stochastic differential equations: the EKF, the UKF, the EnKF, and a PF. The state estimation methods are implemented for non-stiff systems in Matlab, and a numerical experiment is performed for a simulated MFTS. The performance of each state estimation method is evaluated in terms of 1) the computational times for the time- and measurement-updates and 2) the accuracy measured by the MAPE for the state and disturbance estimates.